\documentclass{article}
\usepackage{amsmath,amssymb}

 \def\ivq{\in \! \vee q}

\begin{document}
\title{\bf Some kinds of $(\overline{\in},\overline{\in}\vee
\overline{q})$-fuzzy filters of $BL$-algebras}

\author{ {Xueling Ma}$^a$, {Jianming Zhan}$^{a,*}$, {Wies{\l}aw A. Dudek}$^b$\\ {\small
   $^{a }$ Department of Mathematics,
  Hubei Institute for Nationalities,}\\ {\small Enshi, Hubei Province,
   445000, P. R. China}\\
  {\small  $^b$   Institute of Mathematics and Computer Science,
         Wroc{\l}aw University of Technology, }\\
         {\small   Wybrze\.ze Wyspia\'nskiego 27,
         50-370 Wroc{\l}aw, Poland}\\}

\date{}\maketitle

\begin{flushleft}\rule[0.4cm]{12cm}{0.3pt}
\parbox[b]{12cm}{\small
{\bf Abstract} The concepts  of $(\overline{\in},\overline{\in}
\vee \overline{q})$-fuzzy (implicative, positive implicative  and
fantastic) filters  of $BL$-algebras are introduced and some
related properties are investigated. Some  characterizations  of
these generalized  fuzzy filters  are derived. In particular,  we
describe the relationships  among ordinary fuzzy (implicative,
positive implicative  and fantastic) filters, $(\in,\ivq)$-fuzzy
(implicative, positive implicative and fantastic) filters and
$(\overline{\in},\overline{\in} \vee \overline{q})$-fuzzy
(implicative, positive implicative  and fantastic) filters of
$BL$-algebras. Finally, we prove that a fuzzy set $F$ on a
$BL$-algebra $L$ is an $(\overline{\in},\overline{\in} \vee
\overline{q})$-fuzzy implicative filter of $L$  if and only if it
is both $(\overline{\in},\overline{\in} \vee \overline{q})$-fuzzy
positive  implicative filter and an
$(\overline{\in},\overline{\in} \vee \overline{q})$-fuzzy
fantastic filter. \\ \\

{\it Keywords:} $BL$-algebra; (implicative, positive implicative,
fantastic) filter;
 $(\in,\in\vee q)$-fuzzy  (implicative, positive implicative, fantastic) filter;
   $(\overline{\in},\overline{\in}\vee \overline{q})$-fuzzy
 (implicative, positive implicative, fantastic) filter; Fuzzy
  (implicative, positive implicative, fantastic) filter with thresholds.  \\

  {\it 2000
Mathematics Subject Classification:} 03E72; 03G10; 03B52. }
\rule{12cm}{0.3pt}
\end{flushleft}
\footnote{* Corresponding author. Tel/Fax: 0086-718-8437732.
       \\ \it E-mail addresses:\rm \ \
  zhanjianming@hotmail.com(J. Zhan), dudek@im.pwr.wroc.pl (W.A.
Dudek).}

\subsection*{1. Introduction and Preliminaries}

\paragraph{ }  It is well known that certain information processing, especially
inferences based on certain information, is based on classical
two-valued logic. Due to strict and complete logical foundation
(classical logic), making inference levels. Thus, it is natural and
necessary to attempt to establish some rational logic system as the
logical foundation for uncertain information processing. It is
evident that this kind of logic cannot be two-valued logic itself
but might form a certain extension of two-valued logic. Various
kinds of non-classical logic systems have therefore been extensively
researched in order to construct natural and efficient inference
systems to deal with uncertainty.

Logic appears in a `sacred' form (resp., a `profane') which is
dominant in proof theory (resp., model theory). The role of logic
in mathematics and computer science is twofold-as a tool  for
applications in both areas, and a technique for laying the
foundations. Non-classical logic including many-valued logic,
fuzzy logic, etc., takes the  advantage of the classical logic to
handle information with various facets of uncertainty \cite{36},
such as fuzziness, randomness, and so on. Non-classical logic has
become a formal and useful tool for computer science to deal with
fuzzy information and uncertain information. Fuzziness and
incomparability are two kinds of uncertainties often associated
with human's intelligent activities in the real world, and they
exist not only in the processed object itself, but also in the
course of the object being dealt with.

The  concept of $BL$-algebra was introduced by
H$\mathrm{\acute{a}}$jek's as the algebraic structures for his
Basic Logic \cite{14}. A well known example of a $BL$-algebra is
the interval $[0,1]$ endowed with the structure induced by a
continuous $t$-norm. On the other hand, the $MV$-algebras,
introduced by Chang in 1958 (see \cite{3}), are one of the most
well known classes of $BL$-algebras. In order to investigate the
logic system whose semantic truth-value is given by a lattice, Xu
\cite{32} proposed  the concept of lattice implication algebras
and studied the properties of  filters in such algebras \cite{33}.
Later on, Wang \cite{29} proved that the lattice implication
algebras are categorically equivalent to the $MV$-algebras.
Furthermore, in order to provide an algebraic proof of the
completeness theorem of a formal deductive  system \cite{30}, Wang
\cite{31} introduced the concept of $R_0$-algebras.  In fact, the
$MV$-algebras, G$\ddot{o}$del algebras and product algebras are
the most known classes of $BL$-algebras. $BL$-algebras are further
discussed by many researchers, see
\cite{6,44,15,16,19,26,27,28,39,40,42}. Recent investigations are
concerned with non-commutative generalizations for these
structures. In \cite{12}, Georgescu et al. introduced the concept
of pseudo $MV$-algebras as a non-commutative generalization of
$MV$-algebras. Several researchers discussed the properties of
pseudo $MV$-algebras, see \cite{9,10,18,22,23}. Pseudo
$BL$-algebras are a common extension of $BL$-algebras and pseudo
$MV$-algebras, see \cite{7,11,13,24,38,41}. These structures seem
to be a very general algebraic concept in order to express the
non-commutative reasoning.

After the introduction of fuzzy sets by Zadeh \cite{35}, there
have been a number of generalizations of this fundamental concept.
A new type of fuzzy subgroup, that is, the $(\in,\in\!\vee
q)$-fuzzy subgroup, was introduced in an earlier paper of Bhakat
and Das \cite{1,2} by using the combined notions of
``belongingness" and ``quasicoincidence" of fuzzy points and fuzzy
sets, which was introduced by Pu and Liu \cite{21}. In fact, the
$(\in,\ivq)$-fuzzy subgroup is an important generalization of
Rosenfeld's fuzzy subgroup. It is now natural  to investigate
similar type of generalizations of the existing fuzzy subsystems
with other algebraic structures. With this objective in view,
Davvaz \cite{4} applied this theory to near-rings and obtained
some useful results. Further, Davvaz and Corsini \cite{5}
redefined fuzzy $H_v$-submodule and many valued implications. In
\cite{37}, Zhan et al. also discussed the properties of
interval-valued $(\in,\ivq)$-fuzzy hyperideals in hypernear-rings.
For more details, the reader is referred to
\cite{4,5,17,18,37,38,39}.

In \cite{17} the concepts of $(\in,\in\vee q)$-fuzzy (implicative,
positive implicative  and fantastic) filters in $BL$-algebras are
introduced and related properties are investigated. As a
continuation of this paper, we further discuss this topic in this
paper. In Section 2, we describe the relationships  among ordinary
fuzzy filters, $(\in,\ivq)$-fuzzy filters  and
$(\overline{\in},\overline{\in} \vee \overline{q})$-fuzzy filters
of $BL$-algebras. In Section 3, we divide into three subsections.
In Section 3.1, we describe the relationships  among ordinary
fuzzy  implicative filters, $(\in,\ivq)$-fuzzy  implicative
filters  and $(\overline{\in},\overline{\in} \vee
\overline{q})$-fuzzy implicative filters  of $BL$-algebras.  In
Section 3.2, we describe the relationships  among ordinary fuzzy
positive implicative filters, $(\in,\ivq)$-fuzzy positive
implicative filters  and $(\overline{\in},\overline{\in} \vee
\overline{q})$-fuzzy  positive implicative filters  of
$BL$-algebras.  Further,   the relationships  among ordinary fuzzy
fantastic filters, $(\in,\ivq)$-fuzzy  fantastic filters
and$(\overline{\in},\overline{\in} \vee \overline{q})$-fuzzy
fantastic filters  of $BL$-algebras are considered in Section 3.3.
Finally, in Section 4, we prove that a fuzzy set $F$ of  a
$BL$-algebra $L$  is an $(\overline{\in},\overline{\in} \vee
\overline{q})$-fuzzy implicative filter of $L$  if and only if it
is both $(\overline{\in},\overline{\in} \vee \overline{q})$-fuzzy
positive  implicative filter and an
$(\overline{\in},\overline{\in} \vee \overline{q})$-fuzzy
fantastic filter.

\medskip Recall that an algebra $L=(L,\le,\wedge,\vee,
\odot,\rightarrow, 0,1)$ is a $BL$-algebra if  it is a bounded
lattice such that the following conditions are satisfied:

$(i)$  $(L,\odot,1)$ is  a commutative  monoid,

$(ii)$  $\odot$  and $\rightarrow$  form  an adjoin  pair, i.e.,
$z\le x\rightarrow y $ if and only if  $x\odot  z\le y$ for all
$x,y,z\in L$,

$(iii)$  $x\wedge y=x\odot (x\rightarrow y)$,

$(iv)$  $(x\rightarrow y)\vee (y\rightarrow x)=1$.

\medskip

In any  $BL$-algebra $L$, the following statements are true  (see
\cite{14}):

$(1)$  $x\le y\Leftrightarrow x\rightarrow y=1$,

$(2)$ $x\rightarrow (y\rightarrow z)=(x\odot y)\rightarrow
z=y\rightarrow (x\rightarrow z)$,

$(3)$  $x\odot y\le x\wedge y,$

$(4)$  $x\rightarrow y\le (z\rightarrow x)\rightarrow
(z\rightarrow y),    x\rightarrow y \le (y\rightarrow
z)\rightarrow (x\rightarrow z)$,

$(5)$  $x\rightarrow x'=x''\rightarrow x$,

$(6)$  $x\vee x'=1\Rightarrow x\wedge x'=0$,

$(7)$  $x\vee y=((x\rightarrow y)\rightarrow y)\wedge
((y\rightarrow x)\rightarrow x),$\\
where $x'=x\rightarrow 0$.

\medskip In what follows, $ L$ is  a $BL$-algebra unless
otherwise specified.

\medskip
A non-empty  subset $A$ of  $L$ is called a {\it filter} of $L$ if
$1\in A$ and $\forall x\in A$, $y\in L$, $x\rightarrow y\in A
\Rightarrow y\in A$. It is easy to check that a non-empty subset
$A$ of $L$ is a filter of $L$  if and only if it satisfies: $(i)$
$\forall x,y\in L$,  $x\odot y\in A$; $(ii)$ $\forall x\in A$,
$x\le y\Rightarrow y\in A$ (see \cite{2,27,28}).

A filter $A$ of $L$ is called

{\it implicative } if  $x\rightarrow(z'\rightarrow y)\in A$,
$y\rightarrow z\in A\Rightarrow x\rightarrow z\in A$,

{\it positive implicative } if $x\rightarrow (y\rightarrow z)\in
A$, $x\rightarrow y\in A\Rightarrow x\rightarrow z\in A$,

{\it fantastic } if $z\rightarrow (y\rightarrow  x)\in A$, $z\in
A\Rightarrow ((x\rightarrow y)\rightarrow y)\rightarrow x\in A$\\
(see \cite{15,16,17,26,28,38,39}).

\paragraph{ } We now review some fuzzy logic concepts. A fuzzy set of $L$ is
(see [35]).

\paragraph{Definition 1.1.}  A fuzzy set $F$ of
$L$, i.e., a function $F: L\rightarrow [0,1]$, is called  a {\it
fuzzy filter} of $L$ if for all $x,y\in L$

$(F_1)$  $F(x\odot y)\ge\min\{F(x),F(y)\},$

$(F_2)$ $x\le y\Rightarrow F(x)\le F(y)$.

\paragraph{Definition 1.2.} A fuzzy filter $F$ of $L$ is called

 {\it implicative } if

(F3)  $F(x\rightarrow  z)\ge
\mathrm{\min}\{F(x\rightarrow(z'\rightarrow y)), F(y\rightarrow z)
\}$,  for all $x,y,z\in L$.

(ii)  A fuzzy filter $F$ of $L$ is called a  {\it  fuzzy positive
implicative  filter } of $L$ if it satisfies:

(F4) $F(x\rightarrow z)\ge \mathrm{min}\{F(x\rightarrow(y\rightarrow
z)),F(x\rightarrow y) \}$, for all $x,y,z\in L.$

(iii) A fuzzy filter $F$ of $L$ is called a  {\it  fuzzy fantastic
filter} of $L$ if it satisfies:

(F5)  $F(((x\rightarrow y)\rightarrow y)\rightarrow x)\ge
\mathrm{min}\{F(z\rightarrow(y\rightarrow x)),F(z)\}$, for all
$x,y,z\in L$.

For a fuzzy set $F$ of $L$ and $t\in (0,1]$, the crisp set
$U(F;t)=\{x\in L\,|\,F(x)\ge t\}$
 is called the {\it level subset} of $F$.

\paragraph{Theorem 1.3 [15,16].}  {\it A fuzzy  set
$F$ of  $L$ is a fuzzy   (resp., implicative, positive implicative)
filter of $L$ if and only if
 $U(F;t)(\ne\emptyset)$ is a  (resp., implicative,positive implicative)  filter of $L$ for all $t\in
(0,1]$.}

\paragraph {}By the above Theorem, we can get the following:

\paragraph{Theorem 1.4. }  {\it A fuzzy  set
$F$ of  $L$ is a fuzzy  fantastic  filter of $L$ if and only if
 $U(F;t)(\ne\emptyset)$ is a  fantastic   filter of $L$ for all $t\in
(0,1]$.}

\paragraph{ }A fuzzy set $F$ of a  $BL$-algebra $L$ having the form

$$F(y)=\left\{\begin{array}{ll} t(\ne 0) & \mbox{\ \ \ \ \
if\ \ \ \ \ } y=x\\ 0 & \mbox{\ \ \ \ \ if\ \ \ \ \ } y\ne x
\end{array}\right.$$

is said to be {\it fuzzy point with support $x$ and value $t$ and is
denoted by $U(x;t)$}. A fuzzy point $U(x;t)$ is said to {\it belong
to }(resp. be {\it quasi-coincident with} ) a fuzzy set $F$, written
as $U(x;t)\in F$ (resp. $U(x;t) q F$) if $F(x)\ge t$ (resp.
$F(x)+t>1)$. If $U(x;t)\in F$ or (resp. and ) $U(x;t)q F$, then we
write $U(x;t)\in\vee q$(resp. $\in\wedge q$) $F$. The symbol
$\overline{\in\vee q}$ means that $\in\vee q$ does not hold. Using
the notion of `` membership ($\in$)" and ``quasi-coincidence ($q$)"
of fuzzy points with fuzzy subsets, we obtain the concept of
$(\alpha,\beta)$-fuzzy subsemigroup, where $\alpha$ and $\beta$ are
any two of $\{\in,q, \in\vee q, \in\wedge q\}$ with $\alpha\ne
\in\wedge q$, was introduced in [2]. It is noteworthy that the most
viable generalization of Rosenfeld's fuzzy subgroup is the notion of
$(\in,\in\vee q)$-fuzzy subgroup. The detailed study with
$(\in,\in\vee q)$-fuzzy subgroup has been considered in [1].

\paragraph{Definition 1.5 [17].}  A   fuzzy  set $F$ of $L$ is said to be an
{\it   $(\in,\ivq)$-fuzzy  filter   } of $L$ if for all $t,r\in
(0,1]$ and $x,y\in L$,

 (F6) $U(x; t )\in F$ and $U(y; r )\in F$
 imply $U(x\odot y$;   $\mathrm{ min}\{t ,r\})\in\vee q F$,

 (F7)  $U(x;r)\in F$ implies $U(y;r)\in\vee q
 F$ with $x\le y$.

\paragraph{Definition 1.6 [17].} (i)  An   $(\in,\in\vee
q)$-fuzzy  filter $F$ of $L$  is called
 an  {\it    $(\in,\in\vee
q)$-fuzzy  implicative  filter }  of $L$  if it satisfies:

(F8)  $F(x\rightarrow  z)\ge
\mathrm{\min}\{F(x\rightarrow(z'\rightarrow y)), F(y\rightarrow
z),0.5 \}$,  for all $x,y,z\in L$.

(ii)  An  $(\in,\in\vee q)$-fuzzy filter $F$ of $L$ is called an
{\it   $(\in,\in\vee q)$-fuzzy positive implicative  filter } of $L$
if it satisfies:

(F9) $F(x\rightarrow z)\ge \mathrm{min}\{F(x\rightarrow(y\rightarrow
z)),F(x\rightarrow y),0.5 \}$, for all $x,y,z\in L.$

(iii) An  $(\in,\in\vee q)$-fuzzy filter $F$ of $L$ is called an
{\it   $(\in,\in\vee q)$-fuzzy fantastic filter} of $L$ if it
satisfies:

(F10)  $F(((x\rightarrow y)\rightarrow y)\rightarrow x)\ge
\mathrm{min}\{F(z\rightarrow(y\rightarrow x)),F(z),0.5\}$, for all
$x,y,z\in L$.

\paragraph{Theorem 1.7 [17].}  {\it A fuzzy  set
$F$ of  $L$ is an $(\in,\ivq)$-fuzzy  (resp., implicative, positive
implicative, fantastic)  filter of $L$ if and only if
$U(F;t)(\ne\emptyset)$ is a  (resp., implicative, positive
implicative, fantastic)  filter of $L$ for all $t\in (0,0.5]$.}

\subsection*{2. Generalized fuzzy filters}

\paragraph{ } Consider $J=\{t\vert t\in (0,1]$ and $U(F;t)$ is an empty set or
a filter of $L$\}. we now consider the following questions:
\begin{description}
\item (i) If  $J=(0.5,1]$, what kind of fuzzy   filters of $L$  will be
$F$?
\item(ii) If $J=(\alpha,\beta], (\alpha,\beta\in(0,1])$, whether $F$ will be a kind of fuzzy
  filters of $L$ or not?
\item(iii) Can we give a description for the  relationship between the above generalized fuzzy
  filters ?
\end{description}

\paragraph{Definition 2.1.} A fuzzy set $F$  of $L$
is called an {\it   $(\overline{\in},\overline{\in} \vee
\overline{q})$-fuzzy filter} of $L$ if for all $t,r\in (0,1]$ and
for all $x,y\in L$,

 (F11) $U(x\odot y;\mathrm{ min}\{t,r\})\overline{\in } F$ implies $U(x;t)\overline{\in} \vee
\overline{q} F$ or $U(y;r)\overline{\in} \vee \overline{q} F$,

 (F12) $U(y;r)\overline{\in}
 F$  implies  $U(x;r)\overline{\in} \vee
\overline{q} F$  with $x\le y$.

\paragraph{Example 2.2.} Let $L=\{0,a,b,1\}$,  where $0<a<b<1$.  Then we define
$x\wedge  y=\mathrm{min}\{x,y\},  x\vee y=\mathrm{max}\{x,y\}$,  and
$\odot$ and $\rightarrow$ as follows:

\begin{center}
\begin{tabular}{c|cccc}
          $\odot$ & 0 & $a$ & $b$ & 1  \\ \hline
          0 & 0 & 0 & 0 & 0  \\
          $a$ & 0 & 0 & $a$ & $a$ \\
          $b$ & 0 & $a$ & $b$ & $b$  \\
          1 & 0 & $a$ & $b$ & 1
   \end{tabular}
\ \ \ \ \ \ \ \ \
\begin{tabular}{c|cccc}
          $\rightarrow$ & 0 & $a$ & $b$ & 1  \\ \hline
          0 & 1 & 1 & 1 & 1  \\
          $a$ & $a$ & 1 & 1 & 1  \\
          $b$ & 0 & $a$ & 1 & 1  \\
          1 & 0 & $a$ & $b$ & 1
\end{tabular}
  \end{center}

 It is clear that  $(L,\wedge,\vee ,\odot, \rightarrow, 1)$ is now a $BL$-algebra.
 Define a  fuzzy set  $F$ in $L$ by  $F(0)=0.2,F(a)=0.5$ and $F(1)=F(b)=0.6$.  It is routine to verify that $F$
  is an $(\overline{\in},\overline{\in} \vee
\overline{q})$-fuzzy filter of $L$, but it could neither be a fuzzy
filter of $L$, nor an
 $(\in,\in\vee q)$-fuzzy filter of $L$.

\paragraph{Theorem 2.3.} {\it A fuzzy set $F$  of $L$   is  an
    $(\overline{\in},\overline{\in}\vee \overline{q})$-fuzzy filter   of $L$
if and only if  for any $x,y\in L$,}

(F13) {\it $\mathrm{  max}\{F(x\odot y), 0.5\}\ge \mathrm{
min}\{F(x),F(y)\}$},

(F14) {\it  $\mathrm{  max}\{F(y),0.5\}\ge
 F(x)$ with $x\le y$}.

\paragraph{Proof.}  (F11)$\Rightarrow$ (F13) If there exists $x,y\in L$
such that $\max\{F(x\odot y),0.5\}<t=\min\{F(x),F(y)\}$, then
$0.5<t\le 1, U(x\odot y; t)\overline{\in} F$ and $U(x;t)\in F,
U(y;t)\in F$. By (F11), we have $U(x;t) \overline{q} F$ or $U(y;t)
\overline{q} F$. Then, $(t\le F(x)$ and $t+F(x)\le 1)$ or $(t\le
F(y)$ and $t+F(y)\le 1)$. Thus, $t\le 0.5$, contradiction.

(F13)$\Rightarrow$ (F11) Let $U(x\odot y;\mathrm{
min}\{t,r\})\overline{\in } F$ , then $F(x\odot y)< \min\{t,r\}$.

(a) If $F(x\odot y)\ge \min\{F(x), F(y)\}$, then $\min\{F(x),F(y)\}<
\min\{t,r\}$, and consequently, $F(x)<t$ or $F(y)<r$. It follows
that $U(x;t)\overline{\in}  F$ or $U(y;r)\overline{\in}  F$. Thus,
$U(x;t)\overline{\in} \vee \overline{q} F$ or $U(y;r)\overline{\in}
\vee \overline{q} F$.

(b) If $F(x\odot y)<\min\{F(x),F(y)\}$, then by (F13), we have
$0.5\ge \min\{F(x),F(y)\}$. Putting  $U(x;t)\in F$  or, $U(y;r)\in
F$, then $t\le F(x)\le 0.5$ or $r\le F(y)\le 0.5$. It follows that
$U(x;t) \overline{q} F$ or $U(y;r) \overline{q} F$,  and thus,
$U(x;t)\overline{\in} \vee \overline{q} F$ or $U(y;r)\overline{\in}
\vee \overline{q} F$.

(F12) $\Rightarrow$ (F14) Let $x\le y$, if there exists $x,y\in L$
such that $\max\{F(y),0.5\}<t=F(x)$, then $0.5<t\le 1,
U(y;t)\overline{\in} F$ and $U(x;t)\in F$.  Since
$U(y;t)\overline{\in} F$, by (F12), we have $U(x;t)\overline{q} F$.
Then $t\le F(x)$ and $t+F(x)\le 1$, which implies, $t\le 0.5$,
contradiction.

(F14) $\Rightarrow$ (F12)  Let $U(y;t)\overline{\in} F$ with $x\le
y,$ then $F(y)<t$.

(a) If $F(y)\ge F(x)$, then $F(x)<t$, and consequently,
$U(x;t)\overline{\in} F$. Thus, $U(x;r)\overline{\in} \vee
\overline{q} F$.

(b) If $F(y)<F(x)$, then by (F14), we have, $0.5\ge F(x)$. Let
$U(x;t)\in F$, then $t\le F(x)\le 0.5$. It follows that $U(x;t)
\overline{q} F$, and thus, $U(x;r)\overline{\in} \vee \overline{q}
F$.

This completes the proof.\ \ $\Box{}$

\paragraph{Lemma 2.4 [17].} {\it Let $F$ be a  fuzzy  set of $L$. Then
$U(F;t)(\ne\emptyset)$ is a  filter of $L$ for all $0.5<t\le 1$ if
and only if  it satisfies (F13) and (F14). }

\paragraph{Theorem 2.5. } {\it A fuzzy  set $F$ of $L$ is an $(\overline{\in},\overline{\in} \vee
\overline{q})$-fuzzy   filter of $L$ if and only if
$U(F;t)(\ne\emptyset)$ is a filter for all $0.5<t\le 1$.}

\paragraph {Proof.} This Theorem is an immediate consequence of Theorem
2.3 and Lemma 2.4.\ \ $\Box{}$

\paragraph {Remark 2.6.} Let $F$ be a   fuzzy  set of a  $BL$-algebra  $L$
and $J=\{t\vert t\in (0,1]$ and $U(F;t)$ an empty subset or a
  filter of $L$\}.

  (i) If $J=(0,1]$, then $F$ is an ordinary
fuzzy  filter of $L$ (Theorem 1.3);

(ii) If $J=(0,0.5]$, then $F$ is an  $(\in,\ivq)$-fuzzy  filter of
$L$ (Theorem 1.7);

(iii) If $J=(0.5,1]$, then $F$ is an $(\overline{\in},\overline{\in}
\vee \overline{q})$-fuzzy   filter of $L$ (Theorem 2.5).

\paragraph {} We now  extend the above theory.

\paragraph {Definition 2.7.} Given  $ \alpha , \beta\in  (0,1]$ and
 $ \alpha < \beta $, we call a  fuzzy  set $F$ of $L$ a
   {\it  fuzzy   filter with thresholds $( \alpha , \beta]$} of $L$ if for all $x,y\in L$,
 the following conditions are satisfied :

(F15)  $\mathrm{ max}\{F(x\odot y), \alpha \}\ge \mathrm{
min}\{F(x),F(y), \beta \}$,

(F16) $\mathrm{ max}\{F(y), \alpha  \}\ge \mathrm{ min}\{F(x), \beta
\}$ with $x\le y$.

\paragraph {Theorem 2.8.} {\it  A fuzzy  set $F$ of $L$ is a
  fuzzy   filter  with thresholds $( \alpha ,  \beta]$ of $L$ if and only if
$U(F;t)(\ne\emptyset)$ is a  filter of $L$ for all $ \alpha <t\le
\beta $.}

\paragraph {Proof.} The proof is similar to the proof of Lemma 2.4.\ \ $\Box{}$

\paragraph {Remark 2.9.} (1) By Definition 2.5,  we have the following result: if $ F$ is a
 fuzzy  filter with thresholds
$(\alpha,\beta]$ of $L$, then we can conclude that

(i) $F$ is an ordinary fuzzy implicative filter when $\alpha=0,
\beta=1$;

(ii) $F$ is an  $(\in,\ivq)$-fuzzy  filter when $\alpha=0,
\beta=0.5$;

(iii) $F$ is an   $(\overline{\in},\overline{\in} \vee \overline{q
})$-fuzzy
 filter when  $\alpha=0.5,  \beta=1$.

(2) By Definition 2.5, we can define other  fuzzy  filters of $L$,
 same as the  fuzzy   filter  with thresholds $(0.3,0.9]$,
 with thresholds $(0.4,0.6]$ of $L$, etc.

(3) However,  the   fuzzy  filter  with thresholds of $L$ may not be
the usual fuzzy filter, or may not be an $(\in,\ivq)$-fuzzy  filter,
or may not be  an $(\overline{\in},\overline{\in} \vee \overline{q
})$-fuzzy   filter, respectively. These situations can be shown in
the following example:

 \paragraph{Example 2.10.} Consider the $BL$-algebra $L$ as in Example
 2.2.
Define a   fuzzy  set  $F$ of $L$ by  $F(0)= 0.2, F(a)=0.4,
F(b)=0.8$ and $F(1)=0.6$.

Then, we have

$$  U(F;t)=\left\{\begin{array}{l l} \{0,a,b,1\} & \mbox{\ \ \ \ \
if\ \ \   }  0<t\le 0.2,\\ \{1,b,a\} & \mbox{\ \ \ \ \ if\ \ \ }
0.2<t\le 0.4,\\ \{1,b\} & \mbox{\ \ \ \ \ if\ \ \ }
0.4<t\le 0.6,  \\ \{b\} & \mbox{\ \ \ \ \ if\ \ \ } 0.6<t\le 0.8,\\
\emptyset & \mbox{\ \ \ \ \ if\ \ \ } 0.8<t\le 1.
\end{array}\right.$$

Thus,  $F$ is a  fuzzy    filter with thresholds (0.4, 0.6]  of $L$.
But $F$ could neither be a fuzzy   filter, an $(\in,\ivq)$-fuzzy
  filter of $L$, nor
 an  $(\overline{\in},\overline{\in} \vee \overline{q})$-fuzzy   filter of
 $L$.

\subsection*{3. Generalized fuzzy implicative  (positive
implicative, fantastic)  filters}

\paragraph { } In this  Section, we divide into three
 parts. In Section 3,1, we describe the relationships  among
ordinary fuzzy  implicative filters, $(\in,\ivq)$-fuzzy  implicative
 filters  and $(\overline{\in},\overline{\in} \vee \overline{q})$-fuzzy
 implicative filters  of $BL$-algebras.  In Section 3.2,  we   describe the relationships  among
ordinary fuzzy   positive implicative filters , $(\in,\ivq)$-fuzzy
positive implicative filters  and $(\overline{\in},\overline{\in}
\vee \overline{q})$-fuzzy  positive implicative
 filters  of $BL$-algebras.  Further,   the relationships  among
ordinary fuzzy  fantastic filters , $(\in,\ivq)$-fuzzy  fantastic
 filters  and$(\overline{\in},\overline{\in} \vee \overline{q})$-fuzzy
  fantastic filters  of $BL$-algebras are considered in Section 3.3.

\subsection*{3.1. Generalized fuzzy implicative  filters}

\paragraph{ } Consider $J=\{t\vert t\in (0,1]$ and $U(F;t)$ is an empty set or
an implicative filter of $L$\}. We now consider the following
questions:
\begin{description}
\item (i) If  $J=(0.5,1]$, what kind of fuzzy   implicative  filters of $L$  will be
$F$?
\item(ii) If $J=(\alpha,\beta], (\alpha,\beta\in(0,1])$, whether $F$ will be a kind of fuzzy
  implicative  filters of $L$ or not?
\item(iii) Can we give a description for the  relationship between the above generalized fuzzy
  implicative  filters ?
\end{description}

\paragraph{Definition 3.1.1.} An $(\overline{\in},\overline{\in} \vee
\overline{q})$-fuzzy filter of $L$ is called an {\it
$(\overline{\in},\overline{\in} \vee \overline{q})$-fuzzy
implicative filter} of $L$ if it satisfies:

 (F17)  $\max\{F(x\rightarrow  z),0.5\}\ge
 \min \{F(x\rightarrow(z'\rightarrow y)), F(y\rightarrow z) \}$,  for
all $x,y,z\in L$.

\paragraph{Example 3.1.2.} Let $L=\{0,a,b,1\}$,  where $0<a<b<1$.  Then we define
$x\wedge  y=\mathrm{min}\{x,y\},  x\vee y=\mathrm{max}\{x,y\}$,  and
$\odot$ and $\rightarrow$ as follows:

\begin{center}
\begin{tabular}{c|cccc}
          $\odot$ & 0 & $a$ & $b$ & 1  \\ \hline
          0 & 0 & 0 & 0 & 0  \\
          $a$ & 0 & $a$ & $a$ & $a$ \\
          $b$ & 0 & $a$ & $a$ & $b$  \\
          1 & 0 & $a$ & $b$ & 1
   \end{tabular}
\ \ \ \ \ \ \ \ \
\begin{tabular}{c|cccc}
          $\rightarrow$ & 0 & $a$ & $b$ & 1  \\ \hline
          0 & 1 & 1 & 1 & 1  \\
          $a$ & $0$ & 1 & 1 & 1  \\
          $b$ & 0 & $b$ & 1 & 1  \\
          1 & 0 & $a$ & $b$ & 1
\end{tabular}
  \end{center}

 It is clear that  $(L,\wedge,\vee ,\odot, \rightarrow, 1)$ is now a $BL$-algebra.
 Define a  fuzzy set  $F$ of $L$ by  $F(0)=0.2,F(a)=0.8, F(b)=0$ and $F(1)=0.6$.  It is routine to verify that $F$
  is an $(\overline{\in},\overline{\in} \vee
\overline{q})$-fuzzy implicative filter of $L$, but it could neither
be a fuzzy  implicative filter of $L$, nor an
 $(\in,\in\vee q)$-fuzzy  implicative filter of $L$.

\paragraph{Lemma 3.1.3.} {\it Let $F$ be a  fuzzy  set of $L$. Then
$U(F;t)(\ne\emptyset)$ is an  implicative   filter of $L$ for all
$0.5<t\le 1$ if and only if  it satisfies (F13),(F14) and (F17). }

\paragraph {Proof.}   Assume that $U(F;t ) (\ne\emptyset)$ is an implicative
filter of $L$. Then,it follows from Lemma 2.4  that (F13)  and (F14)
hold.  If there exist $x,y,z\in L$  such that  $
\mathrm{max}\{F(x\rightarrow z), 0.5\}<t=
 \mathrm{min}\{F(x\rightarrow(z'\rightarrow y)), F(y\rightarrow z), 0.5\}
$,  then $0.5< t \le 1 , F(x\rightarrow z)< t $ and
$x\rightarrow(z'\rightarrow y), y\rightarrow z\in U(F;t ).$ Since
$U(F;t )$ is an implicative  filter of $L$, $x\rightarrow z\in U(F;t
)$, and so $F(x\rightarrow z)\ge t $, which is a contradiction.
Hence (F17) holds.

Conversely, suppose that the conditions (F13), (F14) and (F17) hold.
Then, it follows from Lemma 2.4 that $U(F;t )$  is a filter  of $L$.
Assume that  $0.5<t \le   1 , x\rightarrow (z'\rightarrow y),
y\rightarrow z\in U(F;t )$. Then $0.5< t  \le
\mathrm{min}\{F(x\rightarrow(z'\rightarrow y)),F(y\rightarrow z)\}$
$ \le  \mathrm{max}\{F(x\rightarrow z), 0.5\}<F(x\rightarrow z),$
which implies that $x\rightarrow z\in U(F;t )$. Thus  $U(F;t )$ is
an implicative  filter  of $L$.  \ \ $\Box{}$

\paragraph{Theorem 3.1.4. } {\it A fuzzy  set $F$ of $L$ is an $(\overline{\in},\overline{\in} \vee
\overline{q})$-fuzzy implicative  filter of $L$ if and only if
$U(F;t)(\ne\emptyset)$ is an  implicative  filter for all $0.5<t\le
1$.}

\paragraph {Proof.} This Theorem is an immediate consequence of Theorem
2.5 and Lemma 3.1.3.\ \ $\Box{}$

\paragraph {Remark 3.1.5.} Let $F$ be a   fuzzy  set of a  $BL$-algebra  $L$
and $J=\{t\vert t\in (0,1]$ and $U(F;t)$ an empty subset or an
  implicative   filter of $L$\}.

  (i) If $J=(0,1]$, then $F$ is an ordinary
fuzzy  implicative   filter of $L$ (Theorem 1.3);

(ii) If $J=(0,0.5]$, then $F$ is an  $(\in,\ivq)$-fuzzy  implicative
filter of $L$ (Theorem 1.7);

(iii) If $J=(0.5,1]$, then $F$ is an $(\overline{\in},\overline{\in}
\vee \overline{q})$-fuzzy   implicative   filter of $L$ (Theorem
3.1.4).

\paragraph {} We now  extend the above theory.

\paragraph {Definition 3.1.6.} Given  $ \alpha , \beta\in  (0,1]$ and
 $ \alpha < \beta $, we call a  fuzzy  set $F$ of $L$ a
   {\it  fuzzy  implicative  filter with thresholds $( \alpha , \beta]$} of $L$ if
   it satisfies (F15),(F16) and

(F18)  $\max\{F(x\rightarrow  z),\alpha\}\ge
 \min \{F(x\rightarrow(z'\rightarrow y)), F(y\rightarrow z),\beta \}$,  for
all $x,y,z\in L$.

\paragraph {Theorem 3.1.7.} {\it  A fuzzy  set $F$ of $L$ is a
  fuzzy    implicative   filter  with thresholds $( \alpha ,  \beta]$ of $L$ if and only if
$U(F;t)(\ne\emptyset)$ is an   implicative    filter of $L$ for all
$ \alpha <t\le \beta $.}

\paragraph {Proof.} The proof is similar to the proof of Theorem 3.1.4.\ \ $\Box{}$

\paragraph {Remark 3.1.8.} (1) By Definition 3.1.6,  we have the following result: if $ F$ is a
 fuzzy    implicative  filter with thresholds
$(\alpha,\beta]$ of $L$, then we can conclude that

(i) $F$ is an ordinary fuzzy implicative filter when $\alpha=0,
\beta=1$;

(ii) $F$ is an  $(\in,\ivq)$-fuzzy   implicative   filter when
$\alpha=0, \beta=0.5$;

(iii) $F$ is an   $(\overline{\in},\overline{\in} \vee \overline{q
})$-fuzzy  implicative
 filter when  $\alpha=0.5,  \beta=1$.

(2) By Definition 3.1.6, we can define other  fuzzy   implicative
filters of $L$,
 same as the  fuzzy    implicative   filter  with thresholds $(0.3,0.9]$,
 with thresholds $(0.4,0.6]$ of $L$, etc.

(3) However,  the   fuzzy   implicative   filter  with thresholds of
$L$ may not be the usual fuzzy   implicative  filter, or may not be
an $(\in,\ivq)$-fuzzy   implicative   filter, or may not be  an
$(\overline{\in},\overline{\in} \vee \overline{q })$-fuzzy
implicative    filter, respectively. These situations can be shown
in the following example:

 \paragraph{Example 3.1.9.} Consider $BL$-algebra $L$ as in Example
3.1.2. Define a   fuzzy  set  $F$ of $L$ by  $F(0)= 0.4, F(a)=0.8,
F(b)=0.2$ and $F(1)=0.6$.

Then, we have

$$  U(F;t)=\left\{\begin{array}{l l} \{0,a,b,1\} & \mbox{\ \ \ \ \
if\ \ \   }  0<t\le 0.2,\\ \{1,0,a\} & \mbox{\ \ \ \ \ if\ \ \ }
0.2<t\le 0.4,\\ \{1,a\} & \mbox{\ \ \ \ \ if\ \ \ }
0.4<t\le 0.6,  \\ \{a\} & \mbox{\ \ \ \ \ if\ \ \ } 0.6<t\le 0.8,\\
\emptyset & \mbox{\ \ \ \ \ if\ \ \ } 0.8<t\le 1.
\end{array}\right.$$

Thus,  $F$ is a  fuzzy  implicative    filter with thresholds (0.4,
0.6]  of $L$. But $F$ could neither be a fuzzy  implicative filter,
an $(\in,\ivq)$-fuzzy implicative
  filter of $L$, nor
 an  $(\overline{\in},\overline{\in} \vee \overline{q})$-fuzzy implicative    filter of
 $L$.

\subsection*{3.2. Generalized fuzzy  positive implicative filters}

\paragraph{ } Consider $J=\{t\vert t\in (0,1]$ and $U(F;t)$ is an empty set or
a   positive implicative  filter of $L$\}. We now consider the
following questions:
\begin{description}
\item (i) If  $J=(0.5,1]$, what kind of fuzzy   positive implicative   filters of $L$  will be
$F$?
\item(ii) If $J=(\alpha,\beta], (\alpha,\beta\in(0,1])$, whether $F$ will be a kind of fuzzy
  positive implicative   filters of $L$ or not?
\item(iii) Can we give a description for the  relationship between the above generalized fuzzy
  positive implicative   filters ?
\end{description}

\paragraph{Definition 3.2.1.} An $(\overline{\in},\overline{\in} \vee
\overline{q})$-fuzzy filter of $L$ is  called an {\it
$(\overline{\in},\overline{\in} \vee \overline{q})$-fuzzy positive
implicative  filter} of $L$ if it satisfies:

 (F19)  $\max\{F(x\rightarrow  z),0.5\}\ge
 \min \{F(x\rightarrow(y\rightarrow z)), F(x\rightarrow y) \}$,  for
all $x,y,z\in L$.

\paragraph{Example 3.2.2.} Let $L=\{0,a,b,c,1\}$,  where $0<a<b<c<1$.  Then we define
$x\wedge  y=\mathrm{min}\{x,y\},  x\vee y=\mathrm{max}\{x,y\}$,  and
$\odot$ and $\rightarrow$ as follows:

\begin{center}
\begin{tabular}{c|ccccc}
          $\odot$ & 0 & $a$ & $b$ & $c$ & 1  \\ \hline
          0 & 0 & 0 & 0 & 0 & 0 \\
          $a$ & 0 & $a$ & $a$ & $a$ & $a$ \\
          $b$ & 0 & $a$ & $b$ & $a$ & $b$  \\
          $c$ & 0 & $a$ & $a$ & $c$ & $c$  \\
          1 & 0 & $a$ & $b$ & $c$ & 1
   \end{tabular}
\ \ \ \ \ \ \ \ \
\begin{tabular}{c|ccccc}
          $\rightarrow$ & 0 & $a$ & $b$ & $c$ & 1  \\ \hline
          0 & 1 & 1 & 1 & 1 & 1 \\
          $a$ & $0$ & 1 & 1 & 1 & 1  \\
          $b$ & 0 & $c$ & 1 & $c$ & 1  \\
          $c$ & 0 & $b$ & $b$ & 1 & 1  \\
          1 & 0 & $a$ & $b$ & $c$ & 1
\end{tabular}
  \end{center}

 It is clear that  $(L,\wedge,\vee ,\odot, \rightarrow, 1)$ is now a $BL$-algebra.
 Define a  fuzzy set  $F$ of $L$ by  $F(0)=F(c)=0.2,F(a)=0.4, F(b)=0.6$ and $F(1)=0.8$.
 It is routine to verify that $F$
  is an $(\overline{\in},\overline{\in} \vee
\overline{q})$-fuzzy positive implicative  filter of $L$, but it
could neither be a fuzzy  positive implicative  filter of $L$, nor
an
 $(\in,\in\vee q)$-fuzzy  positive implicative  filter of $L$.

\paragraph{Lemma 3.2.3.} {\it Let $F$ be a  fuzzy  set of $L$. Then
$U(F;t)(\ne\emptyset)$ is a  positive implicative    filter of $L$
for all $0.5<t\le 1$ if and only if  it satisfies (F13),(F14) and
(F19). }

\paragraph {Proof.} It is similar to Lemma 3.1.3.   \ \ $\Box{}$

\paragraph{Theorem 3.2.4. } {\it A fuzzy  set $F$ of $L$ is an $(\overline{\in},\overline{\in} \vee
\overline{q})$-fuzzy positive implicative   filter of $L$ if and
only if $U(F;t)(\ne\emptyset)$ is a   positive implicative   filter
for all $0.5<t\le 1$.}

\paragraph {Proof.} This Theorem is an immediate consequence of Theorem
2.5 and Lemma 3.2.3.\ \ $\Box{}$

\paragraph {Remark 3.2.5.} Let $F$ be a   fuzzy  set of a  $BL$-algebra  $L$
and $J=\{t\vert t\in (0,1]$ and $U(F;t)$ an empty subset or a
  positive implicative    filter of $L$\}.

  (i) If $J=(0,1]$, then $F$ is an ordinary
fuzzy  positive implicative    filter of $L$ (Theorem 1.3);

(ii) If $J=(0,0.5]$, then $F$ is an  $(\in,\ivq)$-fuzzy  positive
implicative filter of $L$ (Theorem 1.7);

(iii) If $J=(0.5,1]$, then $F$ is an $(\overline{\in},\overline{\in}
\vee \overline{q})$-fuzzy   positive implicative    filter of $L$
(Theorem 3.2.4).

\paragraph {} We now  extend the above theory.

\paragraph {Definition 3.2.6.} Given  $ \alpha , \beta\in  (0,1]$ and
 $ \alpha < \beta $, we call a  fuzzy  set $F$ of $L$ a
   {\it  fuzzy  positive implicative   filter with thresholds $( \alpha , \beta]$} of $L$ if
   it satisfies (F15),(F16) and

(F20)  $\max\{F(x\rightarrow  z),\alpha\}\ge
 \min \{F(x\rightarrow(y\rightarrow z)), F(x\rightarrow y),\beta \}$,  for
all $x,y,z\in L$.

\paragraph {Theorem 3.2.7.} {\it  A fuzzy  set $F$ of $L$ is a
  fuzzy    positive implicative    filter  with thresholds $( \alpha ,  \beta]$ of $L$ if and only if
$U(F;t)(\ne\emptyset)$ is a positive implicative     filter of $L$
for all $ \alpha <t\le \beta $.}

\paragraph {Proof.} The proof is similar to the proof of Theorem 3.2.4.\ \ $\Box{}$

\paragraph {Remark 3.2.8.} (1) By Definition 3.2.6,  we have the following result: if $ F$ is a
 fuzzy    positive implicative   filter with thresholds
$(\alpha,\beta]$ of $L$, then we can conclude that

(i) $F$ is an ordinary fuzzy positive implicative  filter when
$\alpha=0, \beta=1$;

(ii) $F$ is an  $(\in,\ivq)$-fuzzy   positive implicative    filter
when $\alpha=0, \beta=0.5$;

(iii) $F$ is an   $(\overline{\in},\overline{\in} \vee \overline{q
})$-fuzzy  positive implicative
 filter when  $\alpha=0.5,  \beta=1$.

(2) By Definition 3.2.6, we can define other  fuzzy   positive
implicative filters of $L$,
 same as the  fuzzy    positive implicative    filter  with thresholds $(0.3,0.9]$,
 with thresholds $(0.4,0.6]$ of $L$, etc.

(3) However,  the   fuzzy   positive implicative    filter  with
thresholds of $L$ may not be the usual fuzzy   positive implicative
filter, or may not be an $(\in,\ivq)$-fuzzy   positive implicative
filter, or may not be  an $(\overline{\in},\overline{\in} \vee
\overline{q })$-fuzzy positive implicative     filter, respectively.
These situations can be shown in the following example:

 \paragraph{Example 3.2.9.} Consider the $BL$-algebra $L$ as in Example
3.2.2. Define a   fuzzy  set  $F$ of $L$ by  $F(0)=F(c)= 0.2,
F(a)=0.4, F(b)=0.8$ and $F(1)=0.6$.

Then, we have

$$  U(F;t)=\left\{\begin{array}{l l} \{0,a,b,c,1\} & \mbox{\ \ \ \ \
if\ \ \   }  0<t\le 0.2,\\ \{1,a,b\} & \mbox{\ \ \ \ \ if\ \ \ }
0.2<t\le 0.4,\\ \{1,b\} & \mbox{\ \ \ \ \ if\ \ \ }
0.4<t\le 0.6,  \\ \{b\} & \mbox{\ \ \ \ \ if\ \ \ } 0.6<t\le 0.8,\\
\emptyset & \mbox{\ \ \ \ \ if\ \ \ } 0.8<t\le 1.
\end{array}\right.$$

Thus,  $F$ is a  fuzzy  positive implicative     filter with
thresholds (0.4, 0.6]  of $L$. But $F$ could neither be a fuzzy
positive implicative  filter, an $(\in,\ivq)$-fuzzy positive
implicative
  filter of $L$, nor
 an  $(\overline{\in},\overline{\in} \vee \overline{q})$-fuzzy positive implicative     filter of
 $L$.

\subsection*{3.3. Generalized fuzzy  fantastic filters}

\paragraph{ } Consider $J=\{t\vert t\in (0,1]$ and $U(F;t)$ is an empty set or
a  fantastic filter of $L$\}. We now consider the following
questions:
\begin{description}
\item (i) If  $J=(0.5,1]$, what kind of fuzzy   fantastic   filters of $L$  will be
$F$?
\item(ii) If $J=(\alpha,\beta], (\alpha,\beta\in(0,1])$, whether $F$ will be a kind of fuzzy
  fantastic   filters of $L$ or not?
\item(iii) Can we give a description for the  relationship between the above generalized fuzzy
  fantastic   filters ?
\end{description}

\paragraph{Definition 3.3.1.} An $(\overline{\in},\overline{\in} \vee
\overline{q})$-fuzzy filter of $L$ is called an {\it
$(\overline{\in},\overline{\in} \vee \overline{q})$-fuzzy  fantastic
filter} of $L$ if it satisfies:

 (F21)  $\max\{F(((x\rightarrow y)\rightarrow y)\rightarrow  x),0.5\}\ge
 \min \{F(z\rightarrow(y\rightarrow x)), F(z) \}$,  for
all $x,y,z\in L$.

\paragraph{Example 3.3.2.} Let $L=\{0,a,b, 1\}$,  where $0<a<b<1$.  Then we define
$x\wedge  y=\mathrm{min}\{x,y\},  x\vee y=\mathrm{max}\{x,y\}$,  and
$\odot$ and $\rightarrow$ as follows:

\begin{center}
\begin{tabular}{c|cccc}
          $\odot$ & 0 & $a$ & $b$ & 1  \\ \hline
          0 & 0 & 0 & 0 & 0  \\
          $a$ & 0 & 0 & 0 & $a$ \\
          $b$ & 0 & 0 & $a$ & $b$  \\
          1 & 0 & $a$ & $b$ & 1
   \end{tabular}
\ \ \ \ \ \ \ \ \
\begin{tabular}{c|cccc}
          $\rightarrow$ & 0 & $a$ & $b$ & 1  \\ \hline
          0 & 1 & 1 & 1 & 1  \\
          $a$ & $b$ & 1 & 1 & 1  \\
          $b$ & $a$ & $b$ & 1 & 1  \\
          1 & 0 & $a$ & $b$ & 1
\end{tabular}
  \end{center}

 It is clear that  $(L,\wedge,\vee ,\odot, \rightarrow, 1)$ is now a $BL$-algebra.
 Define a  fuzzy set  $F$ of $L$ by  $F(a)=0.5, F(b)=F(0)=0.2$ and $F(1)=0.8$.
 It is routine to verify that $F$
  is an $(\overline{\in},\overline{\in} \vee
\overline{q})$-fuzzy fantastic  filter of $L$, but it could neither
be a fuzzy  fantastic  filter of $L$, nor an
 $(\in,\in\vee q)$-fuzzy  fantastic  filter of $L$.

\paragraph{Lemma 3.3.3.} {\it Let $F$ be a  fuzzy  set of $L$. Then
$U(F;t)(\ne\emptyset)$ is a  fantastic    filter of $L$ for all
$0.5<t\le 1$ if and only if  it satisfies (F13),(F14) and (F21). }

\paragraph {Proof.} It is similar to Lemma 3.2.3.   \ \ $\Box{}$

\paragraph{Theorem 3.3.4. } {\it A fuzzy  set $F$ of $L$ is an $(\overline{\in},\overline{\in} \vee
\overline{q})$-fuzzy fantastic   filter of $L$ if and only if
$U(F;t)(\ne\emptyset)$ is a   fantastic   filter for all $0.5<t\le
1$.}

\paragraph {Proof.} This Theorem is an immediate consequence of Theorem
2.5 and Lemma 3.3.3.\ \ $\Box{}$

\paragraph {Remark 3.3.5.} Let $F$ be a   fuzzy  set of a  $BL$-algebra  $L$
and $J=\{t\vert t\in (0,1]$ and $U(F;t)$ an empty subset or a
  fantastic    filter of $L$\}.

  (i) If $J=(0,1]$, then $F$ is an ordinary
fuzzy  fantastic    filter of $L$ (Theorem 1.4);

(ii) If $J=(0,0.5]$, then $F$ is an  $(\in,\ivq)$-fuzzy  positive
implicative filter of $L$ (Theorem 1.7);

(iii) If $J=(0.5,1]$, then $F$ is an $(\overline{\in},\overline{\in}
\vee \overline{q})$-fuzzy   fantastic    filter of $L$ (Theorem
5.4).

\paragraph {} We now  extend the above theory.

\paragraph {Definition 3.3.6.} Given  $ \alpha , \beta\in  (0,1]$ and
 $ \alpha < \beta $, we call a  fuzzy  set $F$ of $L$ a
   {\it  fuzzy  fantastic   filter with thresholds $( \alpha , \beta]$} of $L$ if
   it satisfies (F15),(F16) and

(F22)  $\max\{F(((x\rightarrow y)\rightarrow y)\rightarrow
x),\alpha\}\ge
 \min \{F(z\rightarrow(y\rightarrow x)), F(z),\beta \}$,  for
all $x,y,z\in L$.

\paragraph {Theorem 3.3.7.} {\it  A fuzzy  set $F$ of $L$ is a
  fuzzy    fantastic    filter  with thresholds $( \alpha ,  \beta]$ of $L$ if and only if
$U(F;t)(\ne\emptyset)$ is a fantastic     filter of $L$ for all $
\alpha <t\le \beta $.}

\paragraph {Proof.} The proof is similar to the proof of Theorem 3.3.4.\ \ $\Box{}$

\paragraph {Remark 3.3.8.} (1) By Definition 3.3.6,  we have the following result: if $ F$ is a
 fuzzy    fantastic   filter with thresholds
$(\alpha,\beta]$ of $L$, then we can conclude that

(i) $F$ is an ordinary fuzzy fantastic  filter when $\alpha=0,
\beta=1$;

(ii) $F$ is an  $(\in,\ivq)$-fuzzy   fantastic    filter when
$\alpha=0, \beta=0.5$;

(iii) $F$ is an   $(\overline{\in},\overline{\in} \vee \overline{q
})$-fuzzy  fantastic
 filter when  $\alpha=0.5,  \beta=1$.

(2) By Definition 3.3.6, we can define other  fuzzy   fantastic
filters of $L$,
 same as the  fuzzy    fantastic    filter  with thresholds $(0.3,0.9]$,
 with thresholds $(0.4,0.6]$ of $L$, etc.

(3) However,  the   fuzzy   fantastic    filter  with thresholds of
$L$ may not be the usual fuzzy   fantastic filter, or may not be an
$(\in,\ivq)$-fuzzy   fantastic filter, or may not be  an
$(\overline{\in},\overline{\in} \vee \overline{q })$-fuzzy
fantastic  filter, respectively. These situations can be shown in
the following example:

 \paragraph{Example 3.3.9.} Consider the  $BL$-algebra $L$ as in Example
3.1.2. Define a   fuzzy  set  $F$ of $L$ by  $F(a)=0.8, F(0)=0,
F(b)=0.2$ and $F(1)=0.6$.

Then, we have

$$  U(F;t)=\left\{\begin{array}{l l} \{a,b,1\} & \mbox{\ \ \ \ \
if\ \ \   }  0<t\le 0.2,\\ \{1,a\} & \mbox{\ \ \ \ \ if\ \ \ }
0.2<t\le 0.6,  \\ \{a\} & \mbox{\ \ \ \ \ if\ \ \ } 0.6<t\le 0.8,\\
\emptyset & \mbox{\ \ \ \ \ if\ \ \ } 0.8<t\le 1.
\end{array}\right.$$

Thus,  $F$ is a  fuzzy   fantastic     filter with thresholds (0.2,
0.6]  of $L$. But $F$ could neither be a fuzzy   fantastic  filter,
an $(\in,\ivq)$-fuzzy  fantastic
  filter of $L$, nor
 an  $(\overline{\in},\overline{\in} \vee \overline{q})$-fuzzy  fantastic     filter of
of $L$.

\subsection*{4. Relationships among these generalized fuzzy  filters}

\paragraph{ } In this Section, we discuss the relationships among these generalized fuzzy  filters
of $BL$-algebras and obtain an important result.

\paragraph{Lemma 4.1 [16]. } {\it Every implicative filter of $L$ is
a positive implicative filter.}

\paragraph{Lemma 4.2 [16]. } {\it Let $A$ be a filter of $L$. Then $A$ is  an  implicative filter
of $L$ if and only if $(x\rightarrow y)\rightarrow x\in A\Rightarrow
x\in A$, for all $x,y\in L$.}

\paragraph{ }By the definition of fantastic filters of $L$, we can
immediately get the following:

\paragraph{Lemma 4.3.} {\it Let $A$ be a filter of $L$. Then $A$ is
a fantastic filter if and only if $y\rightarrow x\in A\Rightarrow
((x\rightarrow y)\rightarrow y)\rightarrow x\in A$, for all $x,y\in
L$.}

\paragraph{Lemma 4.4 [16]. } {\it Let $A$ be a filter of $L$. Then $A$ is  a positive  implicative filter
of $L$ if and only if $ x\rightarrow (x\rightarrow y)\in
A\Rightarrow x\rightarrow y\in A$, for all $x,y\in L$.}

\paragraph{Lemma 4.5. } {\it Every implicative filter of $L$ is
a fantastic filter.}

\paragraph{Proof. } Let $A$ be an implicative filter of $L$. For any
$x,y\in L$ be such that $y\rightarrow x\in A$. Since
$x\odot((x\rightarrow y)\rightarrow y)\le x$, and so $x\le
((x\rightarrow y)\rightarrow y)\le x$, which implies,
$(((x\rightarrow y)\rightarrow y)\rightarrow x\le x\rightarrow y$.

Thus, $((((x\rightarrow y)\rightarrow y)\rightarrow x)\rightarrow
y)\rightarrow (((x\rightarrow y)\rightarrow y)\rightarrow x)$

$\ge (x\rightarrow y)\rightarrow (((x\rightarrow y)\rightarrow
y)\rightarrow x)$

$\ge ((x\rightarrow y)\rightarrow y)\rightarrow ((x\rightarrow
y)\rightarrow x)$

$\ge y\rightarrow x$.

By hypothesis, we have
 $((((x\rightarrow y)\rightarrow y)\rightarrow x)\rightarrow
y)\rightarrow (((x\rightarrow y)\rightarrow y)\rightarrow x)\in A.$
It follows from Lemma 4.2 that $((x\rightarrow y)\rightarrow x\in
A$. This proves that $y\rightarrow x\in A\Rightarrow ((x\rightarrow
y)\rightarrow y)\rightarrow x\in A$. Thus, by Lemma 4.3, we know $A$
is a fantastic filter of $L$.\ \ $\Box{}$

\paragraph{Theorem 4.6.} {\it  A  non-empty subset $A$ of $L$ is an
implicative filter of $L$ if and only if it is both a positive
implicative filter and a fantastic filter.}

\paragraph{Proof. }  Necessity: Lemma 4.1 and 4.5.

Sufficiency: Let $x,y\in L$ be such that $(x\rightarrow
y)\rightarrow x\in A$. Since $(x\rightarrow y)\rightarrow x\le
(x\rightarrow y)\rightarrow ((x\rightarrow y)\rightarrow y)$, we
have $((x\rightarrow y)\rightarrow y)\in A$. Since $A$ is  a
positive implicative filter of $L$, by Lemma 4.4, we have

$(x\rightarrow y)\rightarrow y\in A.$  \ \ (*)

Since $(x\rightarrow y)\rightarrow x\le y\rightarrow x$, we have
$y\rightarrow x\in A$. By Lemma 4.3, we have

$((x\rightarrow y)\rightarrow y)\rightarrow x\in A$.  \   \ (**)

By (*) and (**), we have $x\in A$ since $A$ is a filter of $L$.

This proves that $(x\rightarrow y)\rightarrow x\in A\Rightarrow x\in
A$. It follows from Lemma 4.2 that $A$ is an implicative filter of
$L$. \ \ $\Box{}$

\paragraph{Corollary 4.7.} {\it  A  non-empty subset $U(F;t)$ of $L$ is an
implicative filter of $L$ if and only if it is both a positive
implicative filter and a fantastic filter for all $t\in (0.5,1]$.}

\paragraph{ } Finally, we give the relationships among $(\overline{\in},\overline{\in} \vee \overline{q})$-fuzzy
 implicative filters, $(\overline{\in},\overline{\in} \vee
\overline{q})$-fuzzy positive  implicative filters and
$(\overline{\in},\overline{\in} \vee \overline{q})$-fuzzy fantastic
filters of $BL$-algebra.

\paragraph{Theorem 4.8.} {\it A fuzzy set $F$ of    $L$ is an
$(\overline{\in},\overline{\in} \vee \overline{q})$-fuzzy
 implicative filter of $L$  if and only if it is both $(\overline{\in},\overline{\in} \vee
\overline{q})$-fuzzy positive  implicative filter and an
$(\overline{\in},\overline{\in} \vee \overline{q})$-fuzzy fantastic
filter.}

\paragraph{Proof. } Let $F$ be an  $(\overline{\in},\overline{\in} \vee \overline{q})$-fuzzy
 implicative filter of $L$. By Theorem 3.1.4, we know  non-empty subset $U(F;t)$ is an
implicative filter of $L$  for all $t\in (0.5,1]$. By Corollary 4.7,
$U(F;t)$ is both a positive  implicative filter and a  fantastic
filter   of $L$ for all $t\in (0.5,1]$.  It follows from Theorem
3.2.4 and 3.3.4 that $F$ is both an $(\overline{\in},\overline{\in}
\vee \overline{q})$-fuzzy positive implicative filter and an
$(\overline{\in},\overline{\in} \vee \overline{q})$-fuzzy fantastic
filter of $L$.

Conversely, assume that  $F$ is both an
$(\overline{\in},\overline{\in} \vee \overline{q})$-fuzzy positive
implicative filter and an $(\overline{\in},\overline{\in} \vee
\overline{q})$-fuzzy fantastic filter of $L$. By Theorem 3.2.4 and
3.3.4, we know  non-empty subset $U(F;t)$ is both a positive
implicative filter and a fantastic filter   of $L$ for all $t\in
(0.5,1]$. By Corollary 4.7, $U(F;t)$ is an implicative filter of $L$
for all $t\in (0.5,1]$. It follows from Lemma 3.1.4 that $F$ is an
$(\overline{\in},\overline{\in} \vee \overline{q})$-fuzzy
 implicative filter of $L$.
\ \ $\Box{}$

\subsection*{Acknowledgements }

\paragraph{ } The research is partially
supported by the Key Science Foundation of Education Committee of
Hubei Province, China (D200729003, D20082903).

\end{document}